\documentclass[12pt,a4paper,onecolumn]{article}
\usepackage{textcomp}
\usepackage{tikz}
\usepackage{amssymb}
\usepackage{caption}
\usepackage{color}
\usepackage{verbatim}
\usepackage{fancybox}
\usepackage{titlesec} 
\usepackage{mathtools}
\usepackage{amsmath}
\usepackage{esvect}
\usepackage{cite}
\usepackage{setspace}
\usepackage{slashbox}
\usepackage{booktabs}
\usepackage{diagbox}

\usepackage[hidelinks]{hyperref}
\usepackage{mathrsfs}
\usetikzlibrary{positioning,fit,calc}
\usepackage{array}
\usepackage{bm}
\usepackage{multirow}
\usepackage{hhline}
\usepackage[super]{nth}
\usetikzlibrary{positioning,decorations.pathreplacing}

\usepackage{mathrsfs}
\usetikzlibrary{positioning,fit,calc}
\usepackage[a4paper,left=1in,top=1in,right=1in,bottom=1in,nohead]{geometry} 
\usepackage[hidelinks]{hyperref}
\usepackage{array}
\usepackage{bm}
\usepackage[symbol]{footmisc} 
\numberwithin{equation}{section}

\usepackage[utf8]{inputenc}
\usepackage{pgfplots}
\usepackage{pgfplotstable}
\pgfplotsset{compat=1.7}
\usepackage{tikz}
\usepackage{graphicx}
\usepackage{caption}
\usepackage{subcaption}

\def\correspondingauthor{\footnote{Corresponding author. Email: peter.frankl@gmail.com.}}

\tikzset{block/.style={draw,thick,text width=2cm,minimum height=1cm,align=center},
         line/.style={-latex}}
\newcolumntype{P}[1]{>{\centering\arraybackslash}m{#1}} 

\titleformat{\section}[block]{\large\scshape\bfseries}{\thesection.}{1em}{} 
\titleformat{\subsection}[block]{\bfseries}{\thesubsection.}{1em}{} 

\newtheorem{defn}{Definition}[section]
\newtheorem{thm}[defn]{Theorem}
\newtheorem{eg}[defn]{Example}
\newtheorem{ppn}[defn]{Proposition}

\newtheorem{cor}[defn]{Corollary}
\newtheorem{lem}[defn]{Lemma}

\newtheorem{clm}[defn]{Claim}

\begin{document}
\pagenumbering{arabic}
\begin{center}
    \textbf{\Large Analogues of Katona's and}\\
\vspace{0.05 in} 
    \textbf{\Large Milner's Theorems for two families}
\vspace{0.1 in} 
    \\{\large Peter Frankl \correspondingauthor{}
\\Alfr{\'e}d R{\'e}nyi Institute of Mathematics,
\\Budapest, Hungary
\vspace{0.1 in} 
\\Willie Wong H.W.
\\National Institute of Education\\Nanyang Technological University, Singapore}
\end{center}

\begin{abstract}
Let $n>s>0$ be integers, $X$ an $n$-element set and $\mathscr{A}, \mathscr{B}\subset 2^X$ two families. If $|A\cup B|\le s$ for all $A\in\mathscr{A}, B\in \mathscr{B}$, then $\mathscr{A}$ and $\mathscr{B}$ are called \textit{cross $s$-union}. Assuming that neither $\mathscr{A}$ nor $\mathscr{B}$ is empty, we prove several best possible bounds. In particular, we show that $|\mathscr{A}|+|\mathscr{B}|\le 1+\sum\limits_{0\le i\le s}{{n}\choose{i}}$. Supposing $n\ge 2s$ and $\mathscr{A},\mathscr{B}$ are antichains, we show that $|\mathscr{A}|+|\mathscr{B}|\le {{n}\choose{1}}+{{n}\choose{s-1}}$ unless $\mathscr{A}=\{\emptyset\}$ or $\mathscr{B}=\{\emptyset\}$. An analogous result for three families is established as well.
\end{abstract}
\section{Introduction}
\indent\par Let $n>s> 0$ be integers and let $X$ be an $n$-element set. Let us use the standard notation: $2^X$ is the power set and ${{X}\choose{i}}$ is the family of all $i$-subsets of $X$. A family $\mathscr{A}\subset 2^X$ is called an \textit{antichain} if $A\subset A'$ never holds for distinct $A,A' \in\mathscr{A}$. The family $\mathscr{A}$ is called \textit{$s$-union} if $|A\cup A'|\le s$ for all $A,A' \in\mathscr{A}$.
\indent\par Analogously, if $\mathscr{A}, \mathscr{B}\subset 2^X$ are two families satisfying $|A\cup B|\le s$ for all $A\in\mathscr{A}, B\in \mathscr{B}$ then they are called \textit{cross $s$-union}. The related notion of \textit{cross $t$-intersecting} describes families satisfying $|A\cap B|\ge t$ for all $A\in\mathscr{A},B\in\mathscr{B}$. Define the \textit{dual} familiy $\mathscr{A}^c=\{X\backslash A: A\in\mathscr{A}\}$. It is easy to see that $\mathscr{A},\mathscr{B}\subset 2^X$ are cross $s$-union if and only if $\mathscr{A}^c$ and $\mathscr{B}^c$ are cross $(n-s$)-intersecting.
\indent\par The notion of $t$-intersecting is encountered more often in the literature, however for the present paper we find it more convenient to deal with cross $s$-union families.
\begin{eg}
Let $a,b$ be non-negative integers satisfying $a+b\le s$. Then the families ${{[n]}\choose{a}}$ and ${{[n]}\choose{b}}$ are cross $s$-union antichains.
\end{eg}
\indent\par We show that this example is the best possible. For a family $\mathscr{F}\subset 2^X$, define the top and bottom sizes $t(\mathscr{F})=\max\{|F|: F\in\mathscr{F}\}$ and $b(\mathscr{F})=\min\{|F|: F\in\mathscr{F}\}$

\begin{thm}\label{thmA3.1.2}
Suppose that $\mathscr{A},\mathscr{B}\subset 2^X$ are cross $s$-union antichains, $n>s>0$. Then, there exist antichains $\tilde{\mathscr{A}},\tilde{\mathscr{B}}\subset 2^X$ so that 
\\(i) $|\mathscr{A}|\le |\tilde{\mathscr{A}}|, |\mathscr{B}|\le |\tilde{\mathscr{B}}|$,
\\(ii) $t(\tilde{\mathscr{A}})\le t(\mathscr{A})$, $t(\tilde{\mathscr{B}})\le t(\mathscr{B})$, and 
\\(iii) $t(\tilde{\mathscr{A}})+t(\tilde{\mathscr{B}})\le s$ 
\\hold. Moreover, unless $\mathscr{A}={{X}\choose{t(\mathscr{A})}}$, $\mathscr{B}={{X}\choose{t(\mathscr{B})}}$, we can choose the above $\tilde{\mathscr{A}},\tilde{\mathscr{B}}$ to satisfy
\\(iv) $|\mathscr{A}|+|\mathscr{B}|<|\tilde{\mathscr{A}}|+|\tilde{\mathscr{B}}|$.
\end{thm}
\begin{cor}\label{corA3.1.3}
Suppose that $\mathscr{A},\mathscr{B}\subset 2^X$ are nonempty cross $s$-union antichains. Then, (i) and (ii) hold.
\\(i)
\begin{align}
|\mathscr{A}|+|\mathscr{B}|\le \max\limits_{0\le i\le \frac{s}{2}}\Bigg\{{{n}\choose{i}}+{{n}\choose{s-i}}\Bigg\}, \label{eqA3.1.1}
\end{align}
with equality if and only if $\{\mathscr{A},\mathscr{B}\}=\{{{X}\choose{i^*}}, {{X}\choose{s-i^*}}\}$ for some integer $0\le i^*\le\frac{s}{2}$.
\\(ii) Assume further $n\ge 2s$ and $\mathscr{A}\neq \{\emptyset\}, B\neq\{\emptyset\}$. Then,
\begin{align}
|\mathscr{A}|+|\mathscr{B}|\le {{n}\choose{1}}+{{n}\choose{s-1}}.\label{eqA3.1.2}
\end{align}
\end{cor}

\indent\par We should mention that (\ref{eqA3.1.1}) is an unpublished result of Ou \cite{OY}. In general it is not easy to determine for which value(s) of $i$ the maximum in (\ref{eqA3.1.1}) is attained. For such $i$, we call $(i,s-i)$ a \textit{maximal pair}. To prove (\ref{eqA3.1.2}), we settle the $n\ge 2s$ case.

\begin{ppn}\label{ppnA3.1.4}
For $n\ge 2s$, the maximal pairs are 
\\(i) $(0,s)$ for all values of $n,s$, except $n=4,s=2$.
\\(ii) $(1,s-1)$ for $n=4,s=2$ and $n=6,s=3$.
\end{ppn}
\indent\par Note that except for the two values in (ii), (\ref{eqA3.1.2}) is a sharpening of (\ref{eqA3.1.1}). We need this stronger bound to establish a result concerning three families.

\begin{thm}\label{thmA3.1.5}
Let $n\ge 2s>0$, with $(n,s)\neq (4,2)$ and suppose that $\mathscr{A},\mathscr{B},\mathscr{C}\subset 2^X$ are nonempty cross $s$-union antichains, that is, $|\mathscr{A}\cup\mathscr{B}\cup\mathscr{C}|\le s$ for all $A\in\mathscr{A}, B\in\mathscr{B}, C\in\mathscr{C}$, then 
\begin{align}
|\mathscr{A}|+|\mathscr{B}|+|\mathscr{C}|\le 2+{{n}\choose{s}}. \label{eqA3.1.3}
\end{align}
with equality if and only if $\{\mathscr{A},\mathscr{B},\mathscr{C}\}=\{{{X}\choose{s}},\{\emptyset\},\{\emptyset\}\}$ unless $(n,s)=(6,3)$, where there is one more optimal triplet $\{\mathscr{A},\mathscr{B},\mathscr{C}\}=\{{{X}\choose{2}},{{X}\choose{1}},\{\emptyset\}\}$.
\end{thm}
We remark that the proof of Theorem \ref{thmA3.1.5} also shows  for $(n,s)=(4,2)$ that $|\mathscr{A}|+|\mathscr{B}|+|\mathscr{C}|\le 9$ and equality holds if and only if $\{\mathscr{A},\mathscr{B},\mathscr{C}\}=\{{{X}\choose{1}},{{X}\choose{1}},\{\emptyset\}\}$.

\indent\par Removing the restriction to be an antichain, we prove another best possible bound.
\begin{thm}\label{thmA3.1.6}
Suppose $\mathscr{A},\mathscr{B}\subset 2^X$ are nonempty and cross $s$-union. Then, for $n>s>0$,
\begin{align} 
|\mathscr{A}|+|\mathscr{B}|\le 1+\sum\limits_{0\le i\le s}{{n}\choose{i}}. \label{eqA3.1.4}
\end{align}
and equality holds if $\{\mathscr{A},\mathscr{B}\}=\{\{\emptyset\},\{F\subset X:\ |F|\le s\}\}$. Furthermore, this pair of extremal families is unique if $s<n-1$.
\end{thm}
\indent\par To put our results in context, let us recall some closely related classical theorems from extremal set theory.
\begin{thm}(Sperner \cite{SE}) 
If $\mathscr{F}\subset2^X$ is an antichain, then $|\mathscr{F}|\le {{n}\choose{\lfloor{n/2}\rfloor}}$ with strict inequality unless $\mathscr{F}={{X}\choose{\lfloor{n/2}\rfloor}}$ or $\mathscr{F}={{X}\choose{\lceil{n/2}\rceil}}$.
\end{thm}

\begin{thm}(Milner \cite{MEC})
Fix $n>s>0$ and let $\mathscr{F}\subset 2^X$ be an $s$-union antichain. Then,
\begin{align}
|\mathscr{F}|\le {{n}\choose{\lfloor {s/2}\rfloor}}. \label{eqA3.1.5}
\end{align}
\end{thm}

\begin{thm}(Frankl \cite{FP 1})\label{thmA3.1.9}
Fix $n>s>0$ and let $\mathscr{A},\mathscr{B}\subset 2^X$ be cross $s$-union antichains. Then,
\begin{align}
\min\{|\mathscr{A}|,|\mathscr{B}|\}\le {{n}\choose{\lfloor{s/2}\rfloor}}. \label{eqA3.1.6}
\end{align}
\end{thm}
Note that setting $\mathscr{A}=\mathscr{B}$, (\ref{eqA3.1.6}) implies (\ref{eqA3.1.5}).

\indent\par There is a recent result which we should state.
\begin{thm}(Wong and Tay \cite{WHW TEG 1})
Suppose that $\mathscr{A},\mathscr{B}\subset 2^X$ are cross $(n-1)$-union antichains. Then,
\begin{align*}
|\mathscr{A}|+|\mathscr{B}|\le {{n}\choose{\lfloor{(n-1)/2}\rfloor}}+{{n}\choose{\lfloor{n/2}\rfloor}},
\end{align*}
with equality if and only if $\{\mathscr{A},\mathscr{B}\}=\big\{{{X}\choose{\lfloor{(n-1)/2}\rfloor}}, {{X}\choose{\lfloor{n/2}\rfloor}}\big\}$.
\end{thm}

\indent\par Recall the following theorem by Katona (see \cite{FP 1} or \cite{OY}) on a $s$-union family.
\begin{thm} (Katona)
If $\mathscr{F}\subset 2^X$ is $s$-union, then $|\mathscr{F}|\le f(n,s)$, where
\begin{align*}
f(n,s)= \left\{
  \begin{array}{@{}ll@{}}
\sum\limits_{0\le i\le r}{{n}\choose{i}}, & \text{if } s=2r \text{ even,} \\[3pt]
2\sum\limits_{0\le i\le r}{{n-1}\choose{i}}, & \text{if } s=2r+1 \text{ odd}. \\
  \end{array}\right.
\end{align*}
\end{thm}

\section{Preliminaries}
\indent\par For a family $\mathscr{F}\subset{{X}\choose{k}}, k\ge 1$, define the \textit{immediate shadow} $\partial\mathscr{F}=\{G\in{{X}\choose{k-1}}:\exists F\in\mathscr{F}, G\subset F\}$; and for $k<n$  the \textit{immediate shade} $\sigma\mathscr{F}=\{H\in{{X}\choose{k+1}}:\exists F\in\mathscr{F}, F\subset H\}$.
\indent\par Recall the following elementary inequalities, due to Sperner \cite{SE}
\begin{align}
&|\partial \mathscr{F}|/|\mathscr{F}|\ge {{n}\choose{k-1}}/{{n}\choose{k}}=\frac{k}{n-k+1},\label{eqA3.2.1}\\
&|\sigma \mathscr{F}|/|\mathscr{F}|\ge {{n}\choose{k+1}}/{{n}\choose{k}}=\frac{n-k}{k+1}.\label{eqA3.2.2}
\end{align}
Furthermore, equality holds in (\ref{eqA3.2.1}) (resp. (\ref{eqA3.2.2})) if and only if $\mathscr{F}=\emptyset$ or  $\mathscr{F}={{X}\choose{k}}$.

\indent\par Define $\mathscr{F}^{(k)}=\{F\in\mathscr{F}: |F|=k\}$. Following Sperner, let us define two new families obtained from $\mathscr{F}$.
\begin{align*}
\mathscr{F}_*=\mathscr{F}\backslash\mathscr{F}^{(t(\mathscr{F}))}\cup\partial(\mathscr{F}^{(t(\mathscr{F}))}),\\
\mathscr{F}^*=\mathscr{F}\backslash\mathscr{F}^{(b(\mathscr{F}))}\cup\sigma(\mathscr{F}^{(b(\mathscr{F}))}).
\end{align*}
\indent\par The next statement is easy to check.
\begin{ppn}\label{ppnA3.2.1}
If $\mathscr{F}$ is a nonempty antichain $(\mathscr{F}\neq\{\emptyset\}$ or $\{X\})$, then both $\mathscr{F}_*$ and $\mathscr{F}^*$ are antichains and $t(\mathscr{F}_*)=t(\mathscr{F})-1$, $b(\mathscr{F}^*)=b(\mathscr{F})+1$.
\end{ppn}
\indent\par The following result is from \cite{FP 1}.
\begin{ppn} (Frankl \cite{FP 1})
Let $k,l$ be positive integers, $\mathscr{G}\subset{{X}\choose{k}}, \mathscr{H}\subset{{X}\choose{l}}$. Suppose that $\mathscr{G}$ and $\mathscr{H}$ are cross $r$-intersecting, $r\ge 1$. Then (i) or (ii) holds.
\\(i) either $|\partial\mathscr{G}|>|\mathscr{G}|$ or $|\partial\mathscr{H}|>|\mathscr{H}|$.
\\(ii) $r=1$, $|\partial\mathscr{G}|=|\mathscr{G}|$, $|\partial\mathscr{H}|=|\mathscr{H}|$.
\end{ppn}
\begin{cor} (Frankl \cite{FP 1}) \label{corA3.2.3}
Suppose that $\mathscr{A},\mathscr{B}\subset 2^X$ are antichains and $\mathscr{A}^{(t(\mathscr{A}))}$ and  $\mathscr{B}^{(t(\mathscr{B}))}$ are cross $r$-intersecting, $r\ge 1$. Then (i) or (ii) holds.
\\(i) either $|\mathscr{A}_*|>|\mathscr{A}|$ or $|\mathscr{B}_*|>|\mathscr{B}|$.
\\(ii) $r=1$, $|\mathscr{A}_*|=|\mathscr{A}|$, $|\mathscr{B}_*|=|\mathscr{B}|$.
\end{cor}

\section{The proofs of Theorem \ref{thmA3.1.2} and its consequences}

\textit{Proof of Theorem \ref{thmA3.1.2}}: Our strategy is simple. We try and replace the pair $(\mathscr{A},\mathscr{B})$ by either $(\mathscr{A}_*, \mathscr{B})$ and $(\mathscr{A}, \mathscr{B}_*)$. In view of Proposition \ref{ppnA3.2.1}, both pairs satisfy (ii) and both are cross $s$-union.
\begin{clm}\label{clmA3.3.1}
We may assume 
\begin{align*}
t(\mathscr{A})+t(\mathscr{B})\le s+1. 
\end{align*}
\end{clm} 
\textit{Proof}: If $t(\mathscr{A})+t(\mathscr{B})\ge s+2$, then for $A\in \mathscr{A}^{(t(\mathscr{A}))}, B\in \mathscr{B}^{(t(\mathscr{B}))}$, we have 
\begin{align*}
|A\cap B|=|A|+|B|-|A\cup B|\ge t(\mathscr{A})+t(\mathscr{B})-s\ge 2.
\end{align*}
\indent\par Now, Corollary \ref{corA3.2.3} implies that for the pair $(\tilde{\mathscr{A}},\tilde{\mathscr{B}})$, $(\tilde{\mathscr{A}},\tilde{\mathscr{B}})=(\mathscr{A}_*, \mathscr{B})$ or $(\tilde{\mathscr{A}},\tilde{\mathscr{B}})=(\mathscr{A}, \mathscr{B}_*)$, (i) and (iv) hold. Thus, we can continue the proof with $(\tilde{\mathscr{A}},\tilde{\mathscr{B}})$ instead of $(\mathscr{A},\mathscr{B})$.
\begin{flushright}
$\Box$
\end{flushright}
\indent\par In view of Corollary \ref{corA3.2.3}(ii), the same argument works in the case $t(\mathscr{A})+t(\mathscr{B})=s+1$, unless $|\mathscr{A}_*|=|\mathscr{A}|$ and $|\mathscr{B}_*|=|\mathscr{B}|$. Assume by symmetry $t(\mathscr{A})\le t(\mathscr{B})$. Consequently,
\begin{align*}
t(\mathscr{A})\le \frac{s+1}{2}\le \frac{n}{2}.
\end{align*}
\indent\par Let us distinguish two cases.
\\Case 1. $\mathscr{A}$ is not uniform, that is, $b(\mathscr{A})< t(\mathscr{A})$. 
\indent\par Now, $t(\mathscr{A}^*)+t(\mathscr{B}_*)=t(\mathscr{A})+t(\mathscr{B})-1=s$, i.e., $\mathscr{A}^*$ and $\mathscr{B}_*$ are cross $s$-union. By (\ref{eqA3.2.2}), $|\mathscr{A}^*|>|\mathscr{A}|$. Since $|\mathscr{B}_*|=|\mathscr{B}|$, we can replace $(\mathscr{A},\mathscr{B})$ by $(\mathscr{A}^*,\mathscr{B}_*)$ and conclude the proof.
\\
\\ Case 2. $\mathscr{A}\subset {{X}\choose{t(\mathscr{A})}}$. 
\indent\par Take an arbitrary set $B_0\in\mathscr{B}^{(t(\mathscr{B}))}$. Fix $A_0\in {{X}\choose{t(\mathscr{A})}}$, $A_0\cap B_0=\emptyset$. Then, $|A_0\cup B_0|=s+1$ implies $A_0\not\in\mathscr{A}$. Consequently, $|\mathscr{A}|<|{{X}\choose{t(\mathscr{A})}}|$. Thus, 
we can replace $(\mathscr{A},\mathscr{B})$ by $\big({{X}\choose{t(\mathscr{A})}},\mathscr{B}_*\big)$ and conclude the proof.
\indent\par In view of the above, we may assume that 
\begin{align*}
t(\mathscr{A})+t(\mathscr{B})\le s.
\end{align*}
\indent\par Suppose again $t(\mathscr{A})\le t(\mathscr{B})$. Thus, $t(\mathscr{A})\le \frac{s}{2}< \frac{n}{2}$.
\indent\par If $\mathscr{A}\neq {{X}\choose{t(\mathscr{A})}}$, then we can replace $(\mathscr{A},\mathscr{B})$ by $\big({{X}\choose{t(\mathscr{A})}},\mathscr{B}\big)$ and we are done. Hence, we may assume $\mathscr{A}={{X}\choose{t(\mathscr{A})}}$.
\indent\par If $t(\mathscr{B})\le \frac{n+1}{2}$, then $|\mathscr{B}|\le |{{X}\choose{t(\mathscr{B})}}|$ follows from (\ref{eqA3.2.2}) and we are done by replacing $\mathscr{B}$ by ${{X}\choose{t(\mathscr{B})}}$.
\indent\par The very final case is $|\mathscr{B}|>\frac{n+1}{2}$. However, then (\ref{eqA3.2.1}) implies $|\mathscr{B}_*|>|\mathscr{B}|$ and we are done by replacing $(\mathscr{A},\mathscr{B})$ by $(\mathscr{A},\mathscr{B}_*)$.
\begin{flushright}
$\Box$
\end{flushright}
\textit{Proof of Corollary \ref{corA3.1.3}}:
\\(i) Let $(\mathscr{A},\mathscr{B})$ be a maximal pair, that is, $\mathscr{A},\mathscr{B}\subset 2^X$ are nonempty cross $s$-union antichains maximizing $|\mathscr{A}|+|\mathscr{B}|$. In view of Theorem \ref{thmA3.1.2}, $\mathscr{A}={{X}\choose{a}}$ and $\mathscr{B}={{X}\choose{b}}$ for some integers $a,b$ satisfying $a+b\le s$.
\indent\par We have to prove that $a+b=s$. Suppose the contrary and assume by symmetry $a\le b$. Then, $a\le \frac{s-1}{2}\le \frac{n}{2}-1$. Consequently, ${{n}\choose{a+1}}>{{n}\choose{a}}$. Replacing $\big({{X}\choose{a}}, {{X}\choose{b}}\big)$ by $\big({{X}\choose{a+1}}, {{X}\choose{b}}\big)$ contradicts maximality.
\\
\\(ii) Let us first prove
\begin{align}
|\mathscr{A}|+|\mathscr{B}|\le \max\limits_{1\le i\le s-1}\Bigg\{{{n}\choose{i}}+{{n}\choose{s-i}}\Bigg\}. \label{eqA3.3.1}
\end{align}
\indent\par In view of Theorem \ref{thmA3.1.2}, the only problem could be, if $t(\tilde{\mathscr{A}})$ or $t(\tilde{\mathscr{B}})$ is exactly $s$,  say $t(\tilde{\mathscr{B}})=s$. If $|\mathscr{A}|\ge 2$, then $|\tilde{\mathscr{A}}|\ge 2$ implies $t(\tilde{\mathscr{A}})\ge 1$ and $t(\tilde{\mathscr{A}})+t(\tilde{\mathscr{B}})\ge s+1$ contradicts Theorem \ref{thmA3.1.2}.
\indent\par The only remaining case is $\mathscr{A}=\{A\}$  for a nonempty set $A\subset X$. Now, the cross $s$-union property implies $A\subset B$ for all $B\in \mathscr{B}^{(s)}$. This in turn implies $|\partial\mathscr{B}^{(s)}|>|\mathscr{B}^{(s)}|$. Replacing $(\mathscr{A},\mathscr{B})$ by $(\mathscr{A},\mathscr{B}_*)$ increases $|\mathscr{A}|+|\mathscr{B}|$ and $t(\mathscr{B}_*)=s-1$. Continuing the proof of Theorem \ref{thmA3.1.2} with $(\mathscr{A},\mathscr{B}_*)$, then leads to (\ref{eqA3.3.1}). Finally, to conclude the proof, we have to prove the inequality (\ref{eqA3.3.6}). We proceed by proving a few other lemmas first.

\begin{lem}\label{lemA3.3.2} For $n\ge 12$,
\\(i) ${{n}\choose{3}}>\frac{n^2}{2}$, and
\\(ii)  ${{n}\choose{5}}>\frac{n^3}{3}$.
\end{lem}
\textit{Proof}: The easy proof is left to the reader.
\begin{flushright}
$\Box$
\end{flushright}

\begin{lem}\label{lemA3.3.3}
For $s\ge 4$, 
\begin{align}
{{2s-1}\choose{1}}+{{2s-1}\choose{s-2}}<{{2s-1}\choose{s-1}}.\label{eqA3.3.2}
\end{align}
\end{lem}
\textit{Proof}: It is easy to verify for $s=4$. We shall prove by induction from $s$ to $s+1$, i.e. we show 
 \begin{align}
{{2s+1}\choose{1}}+{{2s+1}\choose{s-1}}<{{2s+1}\choose{s}}. \label{eqA3.3.3}
\end{align}
Note that
\begin{align*}
2+{{2s-1}\choose{s-3}}&<{{2s-1}\choose{s}}\text{ or equivalently,}\\
{{2s}\choose{0}}+{{2s-1}\choose{0}}+{{2s-1}\choose{s-3}}&<{{2s-1}\choose{s}},
\end{align*} 
and adding to (\ref{eqA3.3.2}), we have
\begin{align*}
{{2s}\choose{0}}+{{2s}\choose{1}}+&{{2s}\choose{s-2}}<{{2s}\choose{s}}.
\end{align*}
Adding ${{2s}\choose{s-1}}$ on both sides, (\ref{eqA3.3.3}) follows.
\begin{flushright}
$\Box$
\end{flushright}

\indent\par For $1\le r\le \frac{n}{2}$, let us introduce the notation $g(n,r)={{n}\choose{r}}-{{n}\choose{r-1}}$.
\begin{lem} \label{lemA3.3.4}
\begin{align}
g(n,r+1)\ &\substack{>\\=\\<}\ g(n,r) \text{ is equivalent to }\label{eqA3.3.4}\\[2pt]
(n-2r)^2\ &\substack{>\\=\\<}\ n+2.\label{eqA3.3.5}
\end{align}
\end{lem}
\textit{Proof}: (\ref{eqA3.3.4}) is equivalent to 
\begin{align*}
{{n}\choose{r+1}}+{{n}\choose{r-1}}\ &\substack{>\\=\\<}\ 2{{n}\choose{r}}.\\
\text{ Dividing by }{{n}\choose{r}}:\ \frac{r}{n-r+1}+\frac{n-r}{r+1}\ &\substack{>\\=\\<}\ 2.
\end{align*}
Multiplying by $(n-r+1)(r+1)$ and rearranging, we have (\ref{eqA3.3.5}) as desired.
\begin{flushright}
$\Box$
\end{flushright}
\indent\par Expressing condition (\ref{eqA3.3.5}) for $r\le \frac{n}{2}$ gives
\begin{align*}
n-\sqrt{n+2}\ \substack{>\\=\\<}\ 2r.
\end{align*}

\begin{lem} For $n\ge 2s$ and $1<i\le\frac{s}{2}$,\label{lemA3.3.5}
\begin{align}
{{n}\choose{i}}+{{n}\choose{s-i}}<{{n}\choose{1}}+{{n}\choose{s-1}},\label{eqA3.3.6}
\end{align}
\end{lem}
\textit{Proof}: It is easy to check that (\ref{eqA3.3.6}) holds for $n\le 11$. Consider $n\ge 12$. For $s<\frac{1}{2}(n-\sqrt{n+2}+4)$, i.e., $2(s-2)<n-\sqrt{n+2}$, using
\begin{align*}
{{n}\choose{i}}-{{n}\choose{1}}=\sum\limits_{t=1}^{i-1}{g(n,t)}\text{ and }{{n}\choose{s-1}}-{{n}\choose{s-i}}=\sum\limits_{t=1}^{i-1}{g(n,s-1-t)},
\end{align*}
(\ref{eqA3.3.6}) follows from the monotonicity of $g(n,r)$ by Lemma \ref{lemA3.3.4}, i.e. $g(n,t)<g(n,s-1-t)$ for $t=1,2,\ldots, i-1$.
\indent\par Suppose $s\ge \frac{1}{2}(n-\sqrt{n+2}+4)$. Then, $s>6$ since $n\ge 12$. Let us try and use induction for moving from $n$ to $n+1$. Since ${{n}\choose{j}}$ is strictly monotone increasing for $0\le j\le \frac{n}{2}$, (\ref{eqA3.3.6}) implies
\begin{align}
{{n}\choose{i}}+{{n}\choose{j}}<{{n}\choose{1}}+{{n}\choose{s-1}} \text{ for all } 1<i<j<s-1,\ i+j\le s. \label{eqA3.3.7}
\end{align}
\indent\par Let us add (\ref{eqA3.3.7}) for the instances 
\begin{align*}
{{n}\choose{i-1}}+{{n}\choose{s-i}}\le{{n}\choose{1}}+{{n}\choose{s-2}}\text{ and } {{n}\choose{i}}+{{n}\choose{s-i-1}}<{{n}\choose{1}}+{{n}\choose{s-2}}.
\end{align*}
Note that the first inequality is not strict as it includes the case $i=2$. So, we have 
\begin{align*}
{{n+1}\choose{i}}+{{n+1}\choose{s-i}}&<2{{n}\choose{1}}+2{{n}\choose{s-2}}\\
&={{n+1}\choose{1}}+{{n+1}\choose{s-1}}+{{n}\choose{1}}-{{n}\choose{0}}-\Big[{{n}\choose{s-1}}-{{n}\choose{s-2}}\Big].
\end{align*}
Since $6<s\le \frac{n}{2}$ and $n\ge 12$, ${{n}\choose{s-2}}>{{n}\choose{5}}>\frac{n^3}{3}$ by Lemma \ref{lemA3.3.2}(ii),
\begin{align*}
{{n}\choose{1}}-{{n}\choose{0}}-\Big[{{n}\choose{s-1}}-{{n}\choose{s-2}}\Big]&=n-1-\frac{n-2s+3}{s-1}{{n}\choose{s-2}}\\
&<n-1-\frac{3}{s-1}\cdot\frac{n^3}{3}\\
&<0.
\end{align*}
So, we have
\begin{align}
{{n+1}\choose{i}}+{{n+1}\choose{s-i}}<{{n+1}\choose{1}}+{{n+1}\choose{s-1}}\label{eqA3.3.8}
\end{align} as desired.
\indent\par There is only one value of $s$ that we did not obtain, namely for $n+1$ even, $s=\frac{n+1}{2}$ (that is $n=2s-1$). We need, for $1<i\le\frac{s}{2}$,
\begin{align}
{{2s}\choose{i}}+{{2s}\choose{s-i}}&<{{2s}\choose{1}}+{{2s}\choose{s-1}}.\label{eqA3.3.9}
\end{align}
\indent\par For $i\neq \frac{s}{2}$, let us add (\ref{eqA3.3.7}) for the instances 
\begin{align*}
{{2s-1}\choose{i-1}}+{{2s-1}\choose{s-i}}&\le{{2s-1}\choose{1}}+{{2s-1}\choose{s-2}}\\
\text{ and } {{2s-1}\choose{i}}+{{2s-1}\choose{s-i-1}}&<{{2s-1}\choose{1}}+{{2s-1}\choose{s-2}}.
\end{align*}
Note that the first inequality is not strict as it includes the case $i=2$. So, we have 
\begin{align*}
{{2s}\choose{i}}+{{2s}\choose{s-i}}&<2{{2s-1}\choose{1}}+2{{2s-1}\choose{s-2}}\\
&={{2s}\choose{1}}+{{2s}\choose{s-1}}+{{2s-1}\choose{1}}-{{2s-1}\choose{0}}-\Big[{{2s-1}\choose{s-1}}-{{2s-1}\choose{s-2}}\Big].
\end{align*}
Since $6<s\le \frac{n}{2}$ and $n\ge 12$, ${{2s-1}\choose{s-2}}>{{n}\choose{5}}>\frac{n^3}{3}$ by Lemma \ref{lemA3.3.2}(ii),
\begin{align*}
{{2s-1}\choose{1}}-{{2s-1}\choose{0}}-\Big[{{2s-1}\choose{s-1}}-{{2s-1}\choose{s-2}}\Big]&=2s-2-\frac{2}{s-1}{{2s-1}\choose{s-2}}\\
&<n-1-\frac{2}{s-1}\cdot\frac{n^3}{3}\\
&<0.
\end{align*}
So, we have (\ref{eqA3.3.9}) as desired.
\indent\par Suppose $i=\frac{s}{2}$. Then, (\ref{eqA3.3.9}) is equivalent to
\begin{align}
2{{2s}\choose{\frac{s}{2}}}&<{{2s}\choose{1}}+{{2s}\choose{s-1}}.\label{eqA3.3.10}
\end{align}
\indent\par Since (\ref{eqA3.3.10}) is easy to verify for $s=4$, it suffices to prove for $s\ge 6$,
\begin{align*} 
2\prod\limits_{0\le i\le \frac{s}{2}}\frac{\frac{s}{2}+1+i}{s+i}<1,
\end{align*}
which follows from
\begin{align*}
2\prod\limits_{0\le i\le \frac{s}{2}}\frac{\frac{s}{2}+1+i}{s+i}<2\Big(\frac{7}{9}\Big)^{\frac{s}{2}+1}<1
\end{align*}
for $s\ge 6$.
\begin{flushright}
$\Box\Box$
\end{flushright}
\indent\par We shall end the section by proving Proposition \ref{ppnA3.1.4}.
\\\textit{Proof of Proposition \ref{ppnA3.1.4}}: It is easy to verify the statement for $n\le 11$. What we need to show for $n\ge 12$ is:
\begin{align*}
{{n}\choose{i}}+{{n}\choose{s-i}}<{{n}\choose{0}}+{{n}\choose{s}} \text{ for }0<i\le \frac{s}{2}.
\end{align*}
In view of Lemma \ref{lemA3.3.5}, it suffices to check 
\begin{align}
{{n}\choose{1}}+{{n}\choose{s-1}}&<{{n}\choose{0}}+{{n}\choose{s}},\label{eqA3.3.11}\\
\text{i.e., }g(n,1)&<g(n,s), \nonumber
\end{align}
which follows from Lemma \ref{lemA3.3.4} if $s<\frac{1}{2}(n-\sqrt{n+2})$.
\indent\par Suppose $s\ge \frac{1}{2}(n-\sqrt{n+2})$. Then, $s>4$ since $n\ge 12$. Hence, ${{n}\choose{s-1}}>{{n}\choose{3}}>\frac{n^2}{2}$  by Lemma \ref{lemA3.3.2}(i). Now,
\begin{align*}
n-1<\frac{1}{s}\cdot\frac{n^2}{2}<\frac{n+1-2s}{s}{{n}\choose{s-1}}
\end{align*}
holds and implies (\ref{eqA3.3.11}).
\begin{flushright}
$\Box$
\end{flushright}

\section{Some inequalities and the proof of Theorem \ref{thmA3.1.5}}
\begin{lem}\label{lemA3.4.1}
Let $n\ge 2s$. Then, (i)-(iii) hold.
\\(i) $3{{s}\choose{\lfloor s/2 \rfloor}}<2+{{n}\choose{s}}$.
\\(ii) $1+{{s}\choose{\lfloor s/2 \rfloor}}+{{n}\choose{s-1}}<2+{{n}\choose{s}}$.
\\(iii) ${{n}\choose{1}}+{{n}\choose{s-2}}+ {{n}\choose{\lfloor (s-1)/2 \rfloor}}<2+{{n}\choose{s}}$ for $s\ge 3$
\end{lem}
\textit{Proof}:
\\(i) Since the LHS is constant with respect to $n$ while the RHS is an increasing function of $n$, it is sufficient to prove the case $n=2s$. For $s=1,2$, one can check the inequality directly. For $s\ge 3$, let us use the identity
\begin{align}
\sum\limits_{0\le i\le s}{{s}\choose{i}}{{s}\choose{s-i}}={{2s}\choose{s}}. \label{eqA3.4.1}
\end{align}
Among the terms on the left, one has ${{s}\choose{\lfloor s/2 \rfloor}}{{s}\choose{\lfloor (s+1)/2 \rfloor}}$. Since ${{s}\choose{\lfloor (s+1)/2 \rfloor}}\ge s\ge 3$, (i) follows.
\\
\\(ii) It is straightforward to check the inequality for $s=1,2,3$. Consider $s\ge 4$. Rearranging (ii), we get
 \begin{align*}
{{s}\choose{\lfloor s/2 \rfloor}}<1+\frac{n+1-2s}{n+1-s}{{n}\choose{s}}.
\end{align*}
Since $\frac{n+1-2s}{n+1-s}$ is decreasing for $0\le s \le \frac{n}{2}$, it suffices to show
\begin{align*}
{{s}\choose{\lfloor s/2 \rfloor}}<\frac{1}{s+1}{{n}\choose{s}}.
\end{align*}
Note that ${{s}\choose{\lfloor (s+1)/2 \rfloor}}\ge s+1$ and (\ref{eqA3.4.1}) does the job again.
\\
\\(iii) Let us settle the case $s=3$ first. Plugging $s=3$ into (iii), we obtain
\begin{align*}
3n<2+{{n}\choose{3}}.
\end{align*}
As for $n\ge 6$, $\frac{{{n}\choose{3}}}{3n}=\frac{(n-1)(n-2)}{18}\ge \frac{5\times 4}{18}>1$, we are done.
\indent\par Let us next show that knowing (iii) for the pairs $(n,s)$ and $(n,s-1)$ implies it for $(n+1,s)$.
\indent\par Indeed, knowing
\begin{align*}
{{n}\choose{1}}+{{n}\choose{\lfloor (s-1)/2\rfloor}}+{{n}\choose{s-2}}<2+{{n}\choose{s}}\text{ and}\\
{{n}\choose{1}}+{{n}\choose{\lfloor (s-2)/2\rfloor}}+{{n}\choose{s-3}}<2+{{n}\choose{s-1}},
\end{align*}
via ${{n}\choose{1}}>{{n}\choose{0}}+2$ and $\frac{s-2}{2}>\frac{s-3}{2}$ implies
\begin{align*}
{{n}\choose{0}}+{{n}\choose{\lfloor (s-1)/2\rfloor-1}}+{{n}\choose{s-3}}<{{n}\choose{s-1}}.
\end{align*}
Now, adding this to the top one gives
\begin{align*}
{{n+1}\choose{1}}+{{n+1}\choose{\lfloor (s-1)/2\rfloor}}+{{n+1}\choose{s-2}}<{{n+1}\choose{s}}
\end{align*}
as desired.
\indent\par Now, we need to prove the base case for the induction, that is, the case $(2s,s)$ with $s\ge 4$. So, 
\begin{align*}
\frac{{{2s}\choose{s-2}}}{{{2s}\choose{s}}}=\frac{(s-1)s}{(s+1)(s+2)}&=1-\frac{4s+2}{(s+1)(s+2)}\text{ and }\\
{{2s}\choose{1}}+{{2s}\choose{\lfloor (s-1)/2\rfloor}}&\le 2{{2s}\choose{\lfloor (s-1)/2\rfloor}}.
\end{align*}
Thus, it is sufficient to show
\begin{align}
{{2s}\choose{s}}\frac{2s+1}{(s+1)(s+2)}\ge {{2s}\choose{\lfloor (s-1)/2\rfloor}}. \label{eqA3.4.2}
\end{align}
\indent\par Let us distinguish two cases according to the parity of $s$.
\\
\\Case 1. $s=2r$.
\indent\par (\ref{eqA3.4.2}) reads now
\begin{align*}
&\frac{(4r)!(4r+1)}{(2r)!(2r+2)!}\ge \frac{(4r)!}{(3r+1)!(r-1)!}.\\
\text{Equivalently, }&(4r+1)\prod\limits_{2r+3\le i\le 3r+1}{i}>\prod\limits_{r\le \tilde{i}\le 2r}{\tilde{i}}=r(r+1)\prod\limits_{r+2\le \tilde{i}\le 2r}{\tilde{i}}.
\end{align*}
Noting that $\frac{3r+1-l}{2r-l}>(\frac{3}{2})$, it is sufficient to show 
\begin{align}
(4r+1)\Big(\frac{3}{2}\Big)^{r-1}>r(r+1).\label{eqA3.4.3}
\end{align}
By Bernoulli's inequality, $(\frac{3}{2})^{r-1}\ge 1+\frac{r-1}{2}=\frac{r+1}{2}$. Now, (\ref{eqA3.4.3}) follows from $\frac{4r+1}{2}>r$.
\\
\\Case 2. $s=2r+1$.
\indent\par It is easy to check that (iii) holds for $n=2s$ when $s=5$. Suppose $s\ge 7$. Again,
\begin{align*}
\frac{{{4r+2}\choose{2r-1}}}{{{4r+2}\choose{2r+1}}}&=1-\frac{8r+6}{2(r+1)(2r+3)}\\
&=1-\frac{4r+3}{(r+1)(2r+3)}\\
&<1-\frac{2}{r+2}.
\end{align*}
Thus, it is sufficient to show 
\begin{align*}
\frac{2}{r+2}{{{4r+2}\choose{2r+1}}}&>2{{{4r+2}\choose{r}}},\text{ or equivalently,}\\
{{{4r+2}\choose{2r+1}}}&>(r+2){{{4r+2}\choose{r}}}.
\end{align*}
Note that
\begin{align*}
{{{4r+2}\choose{2r+1}}}/{{{4r+2}\choose{r}}}=\prod\limits_{0\le i\le r}\frac{2r+2+i}{r+1+i}>\Big(\frac{3r+2}{2r+1}\Big)^{r+1}>\Big(\frac{3}{2}\Big)^{r+1},
\end{align*}
and since $(\frac{3}{2})^{r+1}>r+2$ for all $r\ge 3$, we are done.
\begin{flushright}
$\Box$
\end{flushright}

~\\\textit{Proof of Theorem \ref{thmA3.1.5}:} WLOG, assume $|\mathscr{A}|\ge|\mathscr{B}|\ge |\mathscr{C}|\ge 1$.
\indent\par Should $|\mathscr{C}|=1$ hold, applying (\ref{eqA3.1.1}) and Proposition \ref{ppnA3.1.4} to the cross $s$-union pair $\mathscr{A},\mathscr{B}$ yields 
\begin{align*}
|\mathscr{A}|+|\mathscr{B}|\le 1+{{n}\choose{s}},
\end{align*}
implying (\ref{eqA3.1.3}).
\indent\par From now on, assume $|\mathscr{C}|\ge 2$. Let us distinguish some cases.
\\
\\Case 1. There exist $B\in \mathscr{B}, C\in\mathscr{C}$ such that $|B\cup C|=s$.
\indent\par Then, the cross $s$-union property implies $A\subset B\cup C$ for all $A\in\mathscr{A}$. Thus, $\mathscr{A}\le 
{{s}\choose{\lfloor s/2 \rfloor}}$ yielding for all $s\ge 1$,
\begin{align*}
|\mathscr{A}|+|\mathscr{B}|+|\mathscr{C}|\le 3{{s}\choose{\lfloor s/2 \rfloor}}<2+{{n}\choose{s}},
\end{align*}
where the last inequality follows from Lemma \ref{lemA3.4.1}(i).
\indent\par From now on, we may assume that $\mathscr{B}, \mathscr{C}$ are cross $(s-1)$-union. By Theorem \ref{thmA3.1.9}, 
\begin{align}
|\mathscr{C}|\le {{n}\choose{\lfloor (s-1)/2 \rfloor}}. \label{eqA3.4.4}
\end{align}
Case 2. $\mathscr{A}$ and $\mathscr{B}$ are cross $(s-1)$-union.
\indent\par Now, $|\mathscr{B}|\ge |\mathscr{C}|\ge 2$ and (\ref{eqA3.1.2}) imply $|\mathscr{A}|+|\mathscr{B}|\le {{n}\choose{1}}+{{n}\choose{s-2}}$. Invoking (\ref{eqA3.4.4}),
\begin{align*}
|\mathscr{A}|+|\mathscr{B}|+|\mathscr{C}|\le {{n}\choose{1}}+{{n}\choose{s-2}}+ {{n}\choose{\lfloor (s-1)/2 \rfloor}}.
\end{align*}
Now, (\ref{eqA3.1.3}) follows from Lemma \ref{lemA3.4.1}(iii).
\\
\\Case 3. There exist $A\in\mathscr{A}, B\in\mathscr{B}$ satisfying $|A\cup B|=s$.
\indent\par Consequently, $C\subset A\cup B$ for all $C\in\mathscr{C}$ and this permits to replace (\ref{eqA3.4.4}) by 
\begin{align*}
|\mathscr{C}|\le {{s}\choose{\lfloor s/2 \rfloor}}.
\end{align*}
\indent\par To $\mathscr{A}$ and $\mathscr{B}$ we apply (\ref{eqA3.1.2}):
\begin{align*}
|\mathscr{A}|+|\mathscr{B}|\le 1+{{n}\choose{s-1}}.
\end{align*}
Addiing these inequalities yields
\begin{align*}
|\mathscr{A}|+|\mathscr{B}|+|\mathscr{C}|\le 1+{{s}\choose{\lfloor s/2 \rfloor}}+{{n}\choose{s-1}},
\end{align*}
 which implies (\ref{eqA3.1.3}) with strict inequality by Lemma \ref{lemA3.4.1}(ii). We remark that for $s=1,2$, we have $|\mathscr{C}|\le {{n}\choose{0}}=1$ contrary to our assumptions.
\begin{flushright}
$\Box$
\end{flushright}

\section{The proof of Theorem \ref{thmA3.1.6}}

\indent\par In this section, we are going to use shifting, an important operation invented by Erdos, Ko and Rado \cite{EKR}. For the definition confer \cite{FP 2}, we restrict ourselves to defining shifted families. For convenience we assume $X=\{1,...,n\}$ and let $(a_1,....a_k)$ denote
 the set $\{a_1,...,a_k\}$ when we know that $a_1<...<a_k$. For two sets $A= (a_1,....a_k), B= (b_1,....b_k)$, we say that A precedes B in
 the shifting partial order if $a_i \le b_i$ for all $i$.
\begin{defn}
A family $\mathscr{F}$ is said to be \textbf{shifted} if $F\in \mathscr{F}$ implies $G\in \mathscr{F}$ whenever G precedes F.
\end{defn}
\indent\par The following fact is well-known.
\begin{ppn}
Suppose that $\mathscr{F},\mathscr{G}\subset 2^X$ are cross $s$-union. Then there exist $\mathscr{F}',\mathscr{G}'\subset 2^X$ that are shifted, cross $s$-union and $|\mathscr{F}'|=|\mathscr{F}|, |\mathscr{G}'|=|\mathscr{G}|$.
\end{ppn}
In view of this proposition it is sufficient to prove Theorem \ref{thmA3.1.6} for shifted families.
\\
\\\textit{Proof of Theorem \ref{thmA3.1.6}:} Let us prove (\ref{eqA3.1.4}) using induction on $n$.
\\Base case. $n=s+1$.
\indent\par Obviously, $\mathscr{A}^c$ and $\mathscr{B}$ are disjoint. Thus, $|\mathscr{A}|+|\mathscr{B}|=|\mathscr{A}^c|+|\mathscr{B}|\le 2^{s+1}$, which is the RHS of (\ref{eqA3.1.4}).
\indent\par In the induction step, we employ the shifting technique. WLOG, assume $\mathscr{A},\mathscr{B}\subset 2^X$ are shifted, and $n>s+1$. Set $\mathscr{A}(\bar{n})=\{A\in\mathscr{A}:\ n\not\in A\}$, $\mathscr{A}(n)=\{A\backslash \{n\}: n\in A\in \mathscr{A}\}$ (and do the same for $\mathscr{B}$). Now, $\mathscr{A}(\bar{n}), \mathscr{B}(\bar{n})$ are cross $s$-union antichains on $X\backslash \{n\}$. By the induction hypothesis, 
\begin{align}
|\mathscr{A}(\bar{n})|+|\mathscr{B}(\bar{n})|\le 1+\sum\limits_{0\le i\le s}{{n-1}\choose{i}}.\label{eqA3.5.1}
\end{align}
For $\mathscr{A}(n),\mathscr{B}(n)$, we consider two cases.
\\
\\Case 1. $\mathscr{A}(n)$ or $\mathscr{B}(n)$ is empty. 
\indent\par WLOG, assume $|\mathscr{B}(n)|=0$. Using $|A|\le s-1$ for all $A\in \mathscr{A}(n)$, 
\begin{align*}
|\mathscr{A}(n)|\le \sum\limits_{0\le i\le s-1}{{n-1}\choose{i}}.
\end{align*}
Adding this to (\ref{eqA3.5.1}), $|\mathscr{A}|+|\mathscr{B}|\le 1+\sum\limits_{0\le i\le s}{{n}\choose{i}}$ follows.
\\
\\Case 2. Both $\mathscr{A}(n)$ and $\mathscr{B}(n)$ are nonempty.
\\
\\\textbf{Claim}. $\mathscr{A}(n)$ and $\mathscr{B}(n)$ are cross $(s-2)$-union.
\indent\par Otherwise, we can find $A\in\mathscr{A}(n), B\in \mathscr{B}(n)$, with $|A\cup B|\ge s-1$. If $|A\cup B|\ge s$, then $|(A\cup\{n\}) \cup (B\cup\{n\})|>s$, a contradiction. So far, we have not used the property of shiftedness. Let now $|A\cup B|=s-1$ and pick $x\in [n-1]\backslash (A\cup B)$. By shiftedness, $(A\cup\{x\})\in \mathscr{A}$ and by definition $(B\cup \{n\})\in \mathscr{B}$. However, $|(A\cup\{x\}) \cup (B\cup\{n\})|=|A\cup B|+2=s+1$, a contradiction.
\indent\par Now, the induction hypothesis yields
\begin{align*}
|\mathscr{A}(n)|+|\mathscr{B}(n)|\le 1+\sum\limits_{0\le i\le s-2}{{n-1}\choose{i}}<\sum\limits_{0\le i\le s-1}{{n-1}\choose{i}},
\end{align*}
and (\ref{eqA3.5.1}) follows as in Case 1. However, in this case, we get strict inequality. That is, for $s<n-1$, equality can hold only in Case 1, that is, if either $\mathscr{A}(n)$ or $\mathscr{B}(n)$ is empty. Suppose by symmetry, $\mathscr{B}(n)=\emptyset$. Now $|\mathscr{A}(n)|=\sum\limits_{0\le i\le s-1}{{n}\choose{i}}$ and $|\tilde{A}|<s$ for all $\tilde{A}\in\mathscr{A}(n)$ imply $\mathscr{A}(n)=\{\tilde{A}\subset\{1,2,\ldots,n-1\}:\ |\tilde{A}|<s\}$.
\indent\par By shiftedness, $\tilde{A}\in\mathscr{A}(n)$ and $i\not\in\tilde{A}$ imply $(\tilde{A}\cup\{i\})\in\mathscr{A}$ for $1\le i<n$. We infer $\mathscr{A}=\{A\subset\{1,2,\ldots,n\}:\ |A|\le s\}$. Thus, $\mathscr{B}=\{\emptyset\}$.
\indent\par Since both these families are invariant to shifting, we proved that for $s<n-1$, the pair $\{\{\emptyset\},\{F\subset X:\ |F|\le s\}\}$ is unique to achieve equality in (\ref{eqA3.1.4}) even without assuming shiftedness.
\indent\par In the case $s=n-1$, there are many pairs achieving $|\mathscr{A}|+|\mathscr{B}|=2^n$. Pick a \textit{complex} $\emptyset\neq \mathscr{A}\subsetneqq 2^X$, that is, $A\subset A'\in\mathscr{A}$ implies $A\in\mathscr{A}$. Now, the pair $\{\mathscr{A},2^X\backslash \mathscr{A}^c\}$ is cross $(n-1)$-union with $|\mathscr{A}|+|2^X\backslash \mathscr{A}^c|=2^n$.
\begin{flushright}
$\Box$
\end{flushright}
\textbf{Acknowledgement}
\indent\par The second author would like to thank the National Institute of Education, Nanyang Technological University of Singapore, for the generous support of the Nanyang Technological University Research Scholarship.

\end{document}